\documentstyle[12pt]{article}
\begin{document}

\centerline{\bf  Symplectic Approach of Wess-Zumino-Witten Model}
\centerline{\bf  and Gauge Field Theories }
\vskip.4in
\centerline{\bf Wang Bai-Ling \footnote{Present address: Department of Pure
mathematics, University of Adelaide, Adelaide, SA5005, Australia } }
\vskip.2in
\centerline{Department of Mathematics, Peking University}
\centerline{Beijing 100871, P. R. China}

\vskip .5in
\centerline{Abstract}
\vskip.3in
  {A systematic description of the Wess-Zumino-Witten model is
presented.  The symplectic method plays the major role in this
paper and also gives the relationship between the WZW model and
the Chern-Simons model.  The quantum theory is obtained to give
the projective representation of the Loop group. The Gauss
constraints for the connection whose curvature is only focused
on several fixed points are solved. The Kohno connection and the
Knizhnik-Zamolodchikov equation are derived. The holonomy
representation and $\check R$-matrix representation of braid
group are discussed.}

\newpage
\centerline{\bf Content}
\vskip.3in

  $\S 0.$   Introduction

  $\S 1.$   Classical Phase Space of Wess-Zumino-Witten Model

  $\S 2.$   Symplectic Theory of WZW Model

  $\S 3.$   Geometrical Quantization of Wess-Zumino-Witten Model

  $\S 4.$   Moment Map of Gauge Field Theory

  $\S 5.$   Kohno Connection and Knizhnik-Zamolodchikov Equation

  $\S 6.$   Relations to  $R$-Matrix and Quantum Group

\vskip.2in
\centerline{\bf \S  0. Introduction }
\vskip.2in

   In recent years, the conformal field theories and topological
field theories are the most interesting topics for both
mathematicians and physicists.  The conformal field theories
were initiated by Belavin, Polyakov and Zamolodchikov [4] as
2-dimensional quantum field theory describing the 2-dimensional
critical phenomena. Conformal field theory is characterized by
the symmetries such as Kac-Moody and Virasoro algebra and its
correlation functions are characterized  by differential
equations arising from the representation of infinite
dimensional Lie algebras. A field is called topological if it is
defined for smooth manifolds with no additional structures. As
in Atiyah's paper [1], he axiomatized the topological field
theories  and gave the several outlines to specify the theories.
One aspiring paper may be the Witten's [19], who
interpretated the Jones polynomial of knots in terms of
topological field theories (Chern-Simons theory), where he also
suggested the relations between the Chern-Simons theory and the
Wess-Zumino-Witten theory.

     The Wess-Zumino-Witten model [20] are principal $\sigma $ model with a
topological Wess-Zumino term in the action, which was first introduced
by Witten. There have being tremendous activities inspired by his paper.
I would like to mention Kohno's paper [13], where he derived the Kohno
connection
and its relations with the Knizhink-Zamolodchikov equation, $R$-matrix, and
quantum group. Quantum groups are expected to be the mathematical framework to
describe the symmetry properties of conformal field theory and other
integrable models [5-7,9,15,17,18].

    For conformal field theories, the WZW model is the key point to understand
the rational conformal field theories and their hidden symmetries. The basis
field $g(\xi )$ is a field with values in a compact, semisimple Lie group
$G$, its action is given by (see section 1):
$$
\begin{array}{ll}
S(g)&=-\frac {i}{4\pi }\int _\Sigma Tr(g^{-1}\partial gg^{-1}\overline
\partial g)\\
&\qquad -\frac {i}{12\pi}\int _B Tr(\tilde g^{-1}d\tilde g^{-1} \wedge
\tilde g^{-1}d\tilde g^{-1} \wedge\tilde g^{-1}d\tilde g^{-1})
\end{array}
$$
the first term is the nonlinear $\sigma $-model's action, the second term
is a topological term, which is the integration of the generator of $H^3(G,
\cal Z)$, when $G$ is connected and simple connected. This term is called the
Wess-Zumino term. Through the transgression (see proposition 1 of section 2),
one can obtain a 2-form on a loop space, which is in the same homology class
with symplectic form of the WZW model. It also leads to the equivalence
betweeen
the Wess-Zumino-Witten and the Chern-Simons model.

    The Chern-Simons model is a 2+1 dimensional topological field theory with
the action:
$$
CS(A)=\frac {i}{2\pi }\int _Y Tr(AdA+\frac 23 A^3)
$$
which is the integration of a 3-form over a three dimensional manifold $Y$ and
 equal to the
Yang-Mills action over a four dimensional manifold with boundary $Y$ by the
Stokes
formula, where
$A$ is a $G$-connection of $G$-bundle over $Y$.
 This theory is important because the invariants of manifolds and knots in
three dimensional case can be obtained through the physics calculation.
There is a great surprise, since this 2+1 topological theory has an
intrinsic connection with the 1+1 conformal field theory. In general, one can
sets $\Sigma =\partial Y $ and reduces the connection space of three
 dimensional case to the two dimensional case. Using the Yang-Mills theory on
Riemann surfaces [3], we can get the further results about the gauge field
theory and the conformal
field theory.

     In section 1, we introduce the WZW model and the Chern-Simons model,
derive the symplectic form on the phase space of the WZW model from the
Chern-Simons model. Under certain conditions, the actions of these two models
are equal. For the Chern-Simons model, we use the natural symplectic
form described in M. Atiyah and R. Bott's paper [3].

   In section 2, we investigate the intrinsic geometry of loop group. Under
the transgression isomorphism of proposition 1:
$$
H^3(G, \cal R)\cong H^2(\Omega G, \cal R)
$$

 We know why the  WZW model is equvalent to the Chern-Simons model from
the point of topology. We also show that the Kirillov form is not the
needed one. In finite dimensional case, the natural symplectic on a
co-adjoint orbit is the Kirillov form [12,13], while, for the infinite
dimensional
 case such as the WZW model, the Kirillov form is degenerate. For latter use,
we also quote some results from A. Pressley and G. B. Segal's book [16].
 The main
result in this section is classical description of Kac-Moody and Virasoro
symmetries:
$$
\begin{array}{ll}
\lbrace J_n^a, J_m^b\rbrace & = f_c^{ab} J_{n+m}^c - im\delta _{n, -m}
\delta ^{ab} \\
\lbrace E, J_n^a\rbrace & = -inJ_n^a\\
\lbrace l_n, l_m\rbrace & =
(n - m)l_{n+m}
\end{array}
$$
where the Poisson bracket is defined by the symplectic form $\omega  $ on
$\Omega G$.
 In section 3, we use the geometric quantization to give the
representation of loop group and quantize the WZW model.

 In section 4, we discuss the moment map of gauge field theory under the
three different circumstances: Riemann surfaces without boundary, with
one dimensional  boundary, with n fixed punctures. The third case is most
 interesting, because it leads to the deep relations between the gauge field
theory and the conformal field theory with gauge symmetries. The main
results is to derive the Gauss contraints for the first and third cases
by the Marsden-Weinstein reduction.
 We obtain that, for the latter case, the Gauss contraints are
$$
F^a_{z\overline z}(z, \overline z)=-\frac {4\pi}{k}\sum_{i, a}\delta (z, P_i)
T_a^i
$$
 and also solve this Gauss constraints by using the Green function.

    In section 5, we derive the Kohno connection and prove that this
connection
is flat, its holonomy gives rise to the representation of braid group,
 which is the fundamental group of configuration space,
the parallel displacement gives out the solution of  the Knizhik-Zamolodchikov
equation (KZ equation). This KZ equation was obtained
from Ward identity in [4]
and has a variety of properties which are the fundamental theories for the
conformal theory.

 In the last section, we give a general description about the representation
of braid  group by using the solution of the quantum Yang-Baxter eqation,
$R$-matrix, that is,
$$
R=\sum _{i\neq j}e_{ii}\otimes e_{jj} +q^{\frac 12}\sum _i e_{ii}\otimes
e_{ii}+ (q^{\frac 12}-q^{-\frac 12})\sum _{i < j}e_{ij}\otimes e_{ji}
$$
and show that this representation can be factored through the
Terperley-Lieb-Jones algebra:
$$
\left\{\begin{array}{ll}
&\nu_i\nu_{i\pm 1}\nu_i=\nu_{i\pm 1}\nu_i\nu_{i\pm 1}\\
&\nu_i\nu_j=\nu_j\nu_i \qquad \mid i-j\mid \geq 2\\
&\nu_i^2=(1-q)\nu_i +q
\end{array}\right.
$$
where $q=\exp (-i\frac {2\pi }{k+h})$.
 Due to  Kohno, this representation is equivalent
to the holonomy representation of braid group.

\vskip.2in

\centerline{\bf \S 1. Classical Phase Space of Wess-Zumino-Witten
 Model }
\vskip.2in

    In the 1970's, the non-linear $\sigma$ model
played an important role in the
infinite  conservative flow and the Kac-Moody symmetry of dynamic systems.
 Mathematically, this is equivalent to the minimum surface theory
of  a Riemann surface embedded in a homogeneous space of a certain Lie group,
especially to a Lie group. E. Witten in [20] discussed a  model whose
action is the non-linear $\sigma$ model plus a topological term called
the Wess-Zumino term. From then on, many interesting phenomena occurred in
both mathematics and physics.

    In this section, we will give a general description of the conformal
field theory (Wess-Zumino-Witten
model) and its relations to the 3-dimension gauge theory (Chern-Simons model).
The symplectic structure of the phase space is derived from these relations.

Let $\Sigma$ be a Riemann surface, $G$ a compact Lie group with Lie algebra
$\cal G$ which has a $G$-invariant nondegenerate bilinear form $\langle,
\rangle$ on $\cal G$, sometimes also denoted by $Tr$ for convenience.
The non-linear  $\sigma$ model is defined by the following action:
$$I(g)=-\frac {i}{4\pi} \int_{\Sigma} \langle g^{-1}\partial g , g^{-1}
\overline {\partial} g \rangle$$
where $\partial =\frac {\partial}{\partial z}dz,
{\overline \partial} =\frac {\partial}
{\partial \overline {z}}d\overline {z} $
and  $(z, \overline {z} )$
are the standard complex coordinates on
 $\Sigma $. It is known that there is a right and left bi-invariant, integral
 closed 3-form on $G$:
 $$
 \sigma =\frac {1}{12}Tr(g^{-1}dg \wedge g^{-1}dg \wedge g^{-1}dg)
 $$

 This term has a significant meaning from the view of topology which is
 the generator of $H^{3}(G,\cal {Z} )$
 when $G$ is a connected, simply connected compact
  Lie group. It represents the origin of the topological effects in low
  dimension field theory. The Wess-Zumino-Witten model is given by the
  following action:
$$
S(g)=I(g)-\frac {i}{12\pi  }
\int_{B} Tr(\tilde {g}^{-1}d
\tilde {g}^{-1} \wedge
\tilde {g}^{-1}d\tilde {g}^{-1}
\wedge \tilde {g}^{-1}d\tilde {g}^{-1})$$
where $\tilde {g}$ is an arbitrary extension from $\Sigma$ to  $B$ and $B$ is
a three dimensional manifold whose boundary is $\Sigma$ . Due to the closed
 integral 3-form $\sigma$, a different extension leads to a different action
 by $2k\pi i$ for an integer $k$. So as a complex-valued function on
 the phase space, which consists of all the smooth maps from $\Sigma $
  to $G$, $e^{kS(g)} $ is well-defined. One will see that the
  Euler-Lagrange equation is also not affected by the different extension.

      {\bf Theorem 1}: The Euler-Lagrange equation of classical Wess-
      Zumino-Witten model is given by
   $$\overline {\partial} (\partial g g^{-1}) =0 \Longleftrightarrow \partial
   (g^{-1}\overline {\partial} g)=0$$
  {\bf Proof:}  Let $g_{t}: \Sigma \longrightarrow G$ be a family of
  mappings depending
   on $t$ through a fixed map $g$. The first variation of $S(g_{t})$ is
   given by
   $$
 \begin{array}{ll}
\frac {\partial }{\partial t}\mid _{t=0} S(g_{t})
&=-\frac {i}{4\pi }\int_{\Sigma }\frac {\partial }{\partial t}\mid _{t=0}
\langle g_{t}^{-1}\partial g_{t} , g_{t}^{-1} \overline {\partial }
 g_{t} \rangle \\
&\quad -\frac {i}{12\pi } \int_{B}
 \frac {\partial }{\partial t}\mid _{t=0}
Tr(\tilde {g}_{t}^{-1}d\tilde {g}_{t}^{-1} \wedge
\tilde {g}_{t}^{-1}d\tilde {g}_{t}^{-1}
\wedge \tilde {g}_{t}^{-1}d\tilde {g}_{t}^{-1}) \\
&=-\frac {1}{4\pi }\int_{\Sigma }
\langle -2(\overline {\partial } (\partial g g^{-1})),
(\frac {\partial g_t}{\partial t} g_{t}^{-1}\mid _{t=0}) \rangle
\end{array}
$$
We then  obtained that for any
$\frac {\partial g_{t}}{\partial t} g_{t}^{-1}\vert _{t=0}$,
$\frac {\partial }{\partial t}\vert _{t=0} S(g_{t})=0$ if and only if
$\overline {\partial} (\partial g g^{-1}) =0 $. This is equivalent to
$ \partial (g^{-1}\overline {\partial} g)=0$.

 By Theorem 1, if $g: \Sigma \longrightarrow G $ is a critical point of
 $S(g)$, then $\partial g g^{-1} $ is a holomorphic $\cal G$-valued
 (1,0)-form on $\Sigma $, equivalently, $g^{-1}\overline \partial g$
 is an anti-holomorphic $\cal G $-valued (0,1)-form on $\Sigma$. Because
  the action $S(g)$ is irrelevant to the specific complex structure
  on $\Sigma $, the critical point set  of $S(g)$ is invariant under
  the conformal transformation of $\Sigma $. These lead to the
  conformal symmetry of the Wess-Zumino-Witten model.

    In order to investigate the Kac-Moody symmetry of the WZW model, we
resort to the symplectic techniques in Hamiltonian mechanics. For details,
see [14]. So we need to find out the symplectic structure of the
classical phase space.

    For greater convenience, we take the Riemann surface $\Sigma$ as a disc
$D$
with polar-coordinates $(\rho,\theta)$, i.e. $z=\rho e^{i\theta}$. For a
$\cal C^\infty$-map $g:D\longrightarrow G$, we construct a $\cal G$-valued
1-form on $D$, $A=dgg^{-1}$. The affine space, denoted $\cal A$,
consists of the connection 1-form $A=dgg^{-1}$ derived from a $\cal C^{\infty}$
-map $g:D \longrightarrow G$. We constrict the gauge transformation
with some boundary conditions, that is, for a gauge transformation $h:
D \longrightarrow G$, $h$ must be an identity on the boundary of the disc.
Then the two elements  $dgg^{-1}, dhh^{-1}$ are equivalent under a certian
 gauge
transformation, if and only if $g=h$ on the boundary of $D$. Moreover,
the
connection 1-form is invariant under the constant gauge transformations by
elements of $G$. So the equivalent $\cal G$-valued 1-forms on $D$ is made
up of the phase space which is isomorphic to $LG/G$, where
$LG=\cal C^{\infty }(S^{1},G)$, $S^{1}$ is a unit circle. $LG$ is called
a Loop group which has an outstanding theory described by A.Pressley
and G. Segal
in their book [16].

     From the above, we presume that $LG/G$ is the phase space of
the Wess-Zumino-Witten model. The symplectic form can also be obtained from
the symplectic form of the gauge field theory on the Riemann surface.

     For any $A\in \cal A$, the tangent space of $ \cal A$ at $A$ is still
$\cal A$ because of its affine property. Let $\tilde {\xi },
 \tilde {\eta } \in \cal A$, then  $\tilde {\xi },
 \tilde {\eta }$ are two $\cal G$-valued 1-forms on $D$. The
 natural symplectic form on $\cal A$ is given by
 $$  \omega _{A} ( \tilde {\xi }, \tilde {\eta })
 =\frac {1}{2\pi } \int _{D} Tr(\tilde {\xi }\wedge
  \tilde {\eta }) $$

    This symplectic form is invariant under the gauge transformation.
If we take $\tilde \xi , \tilde \eta$ as the differential
forms derived from two $\cal G$-valued
functions  $\xi , \eta$ on $D$, then
$$\omega _{A} ( \tilde {\xi }, \tilde {\eta })
 =\frac {1}{2\pi } \int _{D} Tr(d\xi \wedge d\eta )=\frac {1}{2\pi }
  \int _{\partial D} Tr(\xi d\eta ) =\frac {1}{2\pi }
  \int _{S^1} <\xi (\theta ), \eta '(\theta )> d\theta$$

     So the reduced symplectic form on $LG/G= \Omega G$, which is
still a group, is a invariant  2-form on $\Omega G$ under the
left-transportation,
whose value at the unit point  is given by
$$\omega _{e} (\xi , \eta )
 =\frac {1}{2\pi} \int _{S^1} <\xi (\theta), \eta '(\theta)>d\theta$$
  where $\xi ,\eta \in T_{e}(\Omega G)
    =\cal C^{\infty} (S^1 ,\cal G )/\cal G $ are two loops through zero
in $\cal G$. As we know, this symplectic form is precisely the
canonical symplectic form on $\Omega G$ as discussed in [16].

    We call $\partial gg^{-1}$ and $g^{-1} \overline \partial g$ as the
left chiral current and the right chiral current respectively. The phase space
which consists of the left chiral current is just
$\Omega G$ with its canonical
symplectic structure as we have shown. For the right chiral current,
the procedure is similar, but the left invariant symplectic structure should
be changed to the right invariant one on the same phase space $\Omega G$.

{\bf Remark:} It is not strange to derive the symplectic structure of
the conformal field theory (Wess-Zumino-Witten model) from the gauge field
theory, partly because they have a natural relation as described above,
and partly because we have an intrinsic formula
between the action of the WZW model
and the action of a special gauge field model (Chern-Simons model).

Let us
say a few words about the  Chern-Simons model. It is a particular example of
the topological gauge field theory, which is based on a gauge field
$A=A^{a}_{\mu}T^{a}dx^{\mu}$ in Lie algebra $\cal G$ with the action
$$CS(A)=\frac {i}{4\pi} \int _{Y}Tr(AdA+\frac{2}{3}A^{3})$$
for a three dimensional manifold Y, where ${T^{a}}$ is  the orthogonal
basis with respect to $G$-invariant bilinear form. $S(A)$ is irrelevant
to the metric of Y. It defines a topological action. Set
$Y=D \times \cal R$, $D$ is a unit disc with polar coordinates $(\rho ,
\theta)$, where $t$ is a coordinate of $\cal R$. Then the d-operator on Y can
be expressed as $\frac {d}{\partial t} dt + \tilde {d}$ with
$\tilde {d}$ the operator on $D$. For g(t): $D
\longrightarrow G $ is a $t$-dependance $C^{\infty}$ map,
$A=-\tilde {d}g(t)g(t)^{-1}$ is a flat connection of the trivial
bundle $(D \times \cal R) \times G$, we have $CS(A)=S_{WZW}(g(t))$.

{\bf Proof:}
$$
\begin{array}{l}
\quad CS(A) \\
=\frac {i}{4\pi }\int_{D\times \cal {R}} Tr(AdA+\frac {2}{3} A^3 )\\
=\frac {i}{4\pi }\int_{D\times \cal {R}}
Tr(A\frac {\partial A}{\partial t}dt) \\
=\frac {i}{4\pi }\int_{D\times \cal {R}}
Tr(\tilde {d}g(t)g(t)^{-1}
\frac {\partial (\tilde {d}g(t)g(t)^{-1})}{\partial t}dt)\\
=\frac {i}{4\pi }\int_{D\times \cal {R}}
Tr(\tilde {d}g(t)g(t)^{-1}\frac {\partial \tilde {d}g(t)}{\partial t}
g(t)^{-1}dt) \\
\qquad -\frac {i}{4\pi }\int_{D\times \cal {R}}
Tr(\tilde {d}g(t)^{-1}\tilde {d}g(t)g(t)^{-1}
\frac {\partial g(t)}{\partial t}g^{-1}(t)dt)\\
=\frac {i}{4\pi }\int_{D\times \cal {R}}
 Tr(-\tilde {d}g(t)^{-1} \tilde {d}
 (\frac {\partial g(t)}{\partial t}dt) -\frac {1}{3}(dgg^{-1})^{3})\\
=-\frac {i}{4\pi }\int_{\partial D \times \cal {R}}
<g(t)^{-1} \frac {\partial g(t)}{\partial \theta }d\theta ,
g(t)^{-1} \frac {\partial g(t)}{\partial t}dt>
-\frac {i}{12\pi }\int_{D\times \cal {R}} Tr(dgg^{-1})^{3}\\
=S_{WZW}(g)
\end{array}
$$
\vskip.3in
\centerline{\bf \S 2. Symplectic Theory of WZW Model }
\vskip.2in
       In this section, we will give a further study  of the phase space
and the symplectic structure from the point of the coadjoint orbit theory.
The left chiral current can be viewed as an infinite function on the
phase space and the conformal symmetry has its role on the phase space.
We also pave the way for a geometrical quantization of the conformal field
theory which  will give rise to representations of the Kac-Moody group and
its Lie algebra, along with the representation of the Virasoro algebra. We
will discuss these in the next section.

      $\Omega G=LG/G$ can be viewed as loops through the unit point in $G$,
its Lie algebra is $L \cal G/ \cal G$, where $L \cal G=C^{ \infty}
(S^{1}, \cal G)$. Another way to view this phase space is that it can
be expressed as a coadjoint orbit of a certain group, which is the
central extension of the Loop group  $LG=\cal C^{\infty }(S^{1}, G)$.

     We define the central extension of the loop algebra $L
\cal G = \cal C^{\infty }(S^{1}, \cal G)$,  as in [16],
by a 2-cocycle $\omega $ which is
exactly $\omega _{e}$ obtained in last section. Specifically, let
$\widetilde {L \cal G} =L \cal G \oplus \cal R$, its Lie bracket is given
by
$$ [(\xi ,\lambda _1), (\eta , \lambda _2)]=([\xi ,\eta ],
\omega (\xi , \eta ))
$$
where $\omega (\xi ,\eta ) =\frac {1}{2\pi } \int _0^{2\pi }
\langle \xi (\theta ), \eta '(\theta )\rangle d\theta$.

     The group corresponding to $\widetilde {L\cal G }$ is the
$U(1)$-extension of the Loop group $LG$, which is $U(1)$ principal bundle
on $LG$  with its first Chern
class $[\omega ]$. For more details, see [16].

     The following proposition gives the relation between the
topological term of WZW action and the symplectic structure on $\Omega G$
which defines the central extension.

{\bf Proposition 1}: Under the transgression isomorphism:
$$
\tau : H^{3}(G, \cal R) \longrightarrow H^{2}(\Omega G, \cal R )
$$
we have $\tau (\sigma )=\omega$ . That is , for $\xi ,\eta \in T_{\gamma}
(\Omega G)$,
$$\tau (\sigma )(\xi ,\eta )=\frac {1}{2\pi }\int _0^{2\pi }\langle
\gamma ^{-1}(\theta )\gamma '(\theta ), [\gamma ^{-1}\xi (\theta ),
\gamma ^{-1}(\theta )\eta (\theta )]\rangle d\theta $$
which is different from $\omega $  by an exact form $d\beta $, where
$\beta $ is 1-form on $\Omega G $ and for $\xi (\theta ) \in T_{\gamma}
(\Omega G)$

$$\beta (\xi )= \frac {1}{4\pi }\int _0^{2\pi }
<\gamma ^{-1}(\theta )\gamma '(\theta ),
\gamma ^{-1}(\theta) \xi (\theta )>d\theta
$$
{\bf Proof:} First, we recall the transgression $\tau $. There is an evaluation
map $\phi : S^1 \times \Omega G \longrightarrow G$ such that a 3-form
$\sigma $ on $G$ can be pulled back to $\Omega G \times S^1 $, then it
can be
integrated over $S^1 $. One then obtains $\tau (\sigma )$.
$$
\phi ^{*}(\sigma )(\frac {\partial }{\partial \theta }, \xi (\theta ),
\eta (\theta )) =\frac {1}{2\pi }\langle
\gamma ^{-1}(\theta )\gamma '(\theta ), [\gamma ^{-1}\xi (\theta ),
\gamma ^{-1}(\theta )\eta (\theta )]\rangle
$$
Integrating it over $S^1$, we obtain $\tau (\sigma )$.

  Second, we need to prove $\omega =\tau (\sigma )-d\beta $. By the
definition of $d \beta  $, we have
$$
\begin{array}{l}
\quad  d\beta _{\gamma }(\xi (\theta ), \eta (\theta ) ) \\
= \xi (\theta ).\beta _{\gamma }(\eta (\theta ))-
\eta (\theta).\beta _{\gamma }(\xi (\theta ))-\beta _{\gamma }
([\xi (\theta ), \eta (\theta)])
\end{array}
$$

  Let $\gamma _t (\theta )$  be a curve through
 $\gamma _0 (\theta ) =\gamma$
with $\frac {\partial }{\partial t } \mid _{t=0} \gamma _t (\theta )
=\xi (\theta )$, so
$$\begin{array}{ll}
\xi (\theta ).\beta _{\gamma }(\eta (\theta ))
&=\xi (\theta ).\int _0^{2\pi}\langle \gamma ^{-1}\gamma ',
\gamma^{-1}\eta (\theta )\rangle d\theta \\
&=\frac {d}{d\theta }\mid _{t=0}\int _0^{2\pi }<\gamma _t^{-1}\gamma '_t,
\gamma _t^{-1}\eta (\theta )>d\theta\\
&=-\int _0^{2\pi } (<\gamma ^{-1} \xi (\theta )\gamma ^{-1}\gamma ',
\gamma ^{-1}\eta (\theta )>+<\gamma ^{-1} \xi ', \gamma ^{-1}\eta >\\
&\qquad -<\gamma ^{-1}\gamma ', \gamma ^{-1}\xi \gamma ^{-1}\eta >)d\theta
\end{array}
$$
Similarly, let $\gamma _t (\theta )$  be a curve through
$\gamma _0 (\theta )=\gamma $
with $\frac {\partial }{\partial t } \mid _{t=0} \gamma _t (\theta)
=\eta (\theta )$, then
$$\begin{array}{ll}
\eta (\theta ).\beta _{\gamma }(\xi (\theta ))
&=\eta (\theta ).\int _0^{2\pi}\langle \gamma ^{-1}\gamma ',
\gamma^{-1}\xi (\theta )\rangle d\theta\\
&=\frac {d}{d\theta }\mid _{t=0}\int _0^{2\pi }<\gamma _t^{-1}\gamma '_t,
\gamma _t^{-1}\xi (\theta )>d\theta\\
&=-\int _0^{2\pi } (<\gamma ^{-1} \eta (\theta )\gamma ^{-1}\gamma ',
\gamma ^{-1}\xi (\theta )>+<\gamma ^{-1} \eta ', \gamma ^{-1}\xi > \\
&\qquad -<\gamma ^{-1}\gamma ', \gamma ^{-1}\eta \gamma ^{-1}\xi >)d\theta
\end{array}
$$
Therefore :
$$
\begin{array}{l}
\quad d\beta _{\gamma }(\xi (\theta ), \eta (\theta ) )  \\
=\frac {1}{4\pi }\int _0^{2\pi }(<\gamma ^{-1}\xi ', \gamma ^{-1}\eta >
-<\gamma ^{-1}\eta ', \gamma ^{-1}\xi > +<\gamma ^{-1}\gamma ',
 [\gamma ^{-1}\xi ,\gamma ^{-1}\eta ]>)d\theta \\
=- \omega _{\gamma }(\xi ,\eta )+\frac {1}{2\pi }\int _0^{2\pi }
<\gamma ^{-1}\gamma ', [\gamma ^{-1}\xi ,\gamma ^{-1}\eta ]>d\theta \\
=-\omega _{\gamma } (\xi ,\eta )+\tau (\sigma )_{\gamma }(\xi ,\eta )
\end{array}
$$
Thus we obtain $\omega (\gamma )=\tau (\sigma ) (\gamma )-
d \beta  (\gamma )$.

{\bf Proposition 2}:  (1) The adjoint action of $LG$ on
$\widetilde {L\cal G}$ is given by
$$Ad\gamma .(\xi , \lambda )= (Ad\gamma .\xi , \lambda -<\gamma ^{-1}
\gamma ', \xi >) $$
where $ Ad \gamma . \xi $ denotes the adjoint action of $ \gamma \in L\cal G $
on $ \xi \in L\cal G $.

(2)The coadjoint action of $ LG $ on $ (\widetilde {L\cal G})^{*} $ is given
by
$$
Ad^{*}\gamma .(\phi ,\lambda) = (Ad^{*}\gamma .\phi
+ \lambda \gamma '\gamma ^{-1}, \lambda)
$$
where $Ad^* \gamma .\phi $ denotes the coadjoint action of $ \gamma \in LG $
on $ \phi \in (L\cal G)^{*} $.

(3)The coadjoint orbit through (0,1), which belongs to $(\widetilde
{L\cal G})^* $, is isomorphic to $\Omega G$.

{\bf Proof:}(1) and (2) are proved in [16]. Therefore, from  (2), we
 know that the
coadjoint
orbit through (0,1) is
$$
\lbrace Ad^*\gamma .(0,1)\mid \gamma \in LG \rbrace
= \lbrace (\gamma '\gamma ,1)
^{-1} )\mid \gamma \in LG \rbrace \cong \Omega G
$$

       In fact, $\Omega G$ is a K$\ddot a$hler manifold stated in the
following theorem. It seems similar to the finite dimensional case which, in
that case,
a regular coadjoint orbit of a compact Lie group is a K$\ddot a$hler
manifold whose K$\ddot a$hler form is precisely the Krillov form, but there
is a slight difference here.
The Krillov form on $\Omega G$ is $\tau (\sigma)$,
which is degenerate, while $\omega _{\gamma }$ is the correct K$\ddot a$hler
form. In  proposition 1, we have given the exact difference.

{\bf Lemma }: ([16])  If we define a left invariant almost complex structure on
$\Omega G $ by $J_g = l_{g*} J_e $, where $J_e $ valued at the
unit point is given
by
$$
J_e:T_e^{\cal C} (\Omega G) \longrightarrow T_e^{\cal C} (\Omega G) = L\cal
G_{\cal C} / \cal G_{\cal C}
$$
$$
J_e(\sum _{k\not= 0}\xi _k z^k) = \sum _{k\not= 0} i sign (k) \xi _k z^k
$$
then $J$ is integral and compatible with $\omega$ in the following meanings:
$$
\begin{array}{ll}
\omega (J\xi ,J\eta ) &= \omega (\xi ,\eta )\\
\omega (\xi ,J\xi)& \geq 0
\end{array}
$$
where $\omega (\xi ,J\xi ) = 0$ if and only if $\xi = 0$.

The proof is straight forward, therefore $(\Omega G, \omega _L ,J_L)$ is a
K$\ddot a$hler manifold.

Let $\lbrace T_a \rbrace _{a = 1}^{dim \cal G}$ be a unitary orthogonal
basis of $\cal G $.  By the Fourier expansion,
$$
\gamma ' (\theta ) \gamma ^{-1}(\theta )
= \sum _{n = - \infty }^{n = +\infty }
J_n^a (\gamma ) T_a e^{-in\theta }
$$
where $J_n^a(\gamma ) \in \cal C^{\infty } (\Omega G )$ is given by
$$
J_n^a (\gamma ) = \frac {1} {2\pi } \int _0 ^{2\pi } <T_ae^{in\theta },
\gamma '(\theta )\gamma ^{-1} (\theta )> d\theta
$$
In addition, set
$$
E(\gamma ) = \frac {1} {2\pi } \int _0 ^{2\pi } < \gamma '(\theta )
\gamma ^{-1} (\theta ),\gamma '(\theta )\gamma ^{-1} (\theta ) > d\theta
$$
$$
l_n(\gamma ) = \frac {1} {4\pi } \int _0 ^{2\pi } < \gamma '(\theta )
\gamma ^{-1} (\theta ) e^{in\theta }, \gamma '(\theta )\gamma ^{-1}
(\theta ) > d\theta
$$

{\bf Theorem 2}: (1) For a smooth function $f$, its Hamiltonian vector field
$X_f $ is given by $X_f\rfloor \omega = - df$, then
$$
X_{J_n^a(\gamma )} = - R_{\gamma *}e^{in\theta }T^a
$$
$$
X_{l_n} (\gamma ) = - e^{in\theta } \gamma '(\theta )
$$
(2)Assume $[T^a, T^b] = f_c^{ab} T^c$. We give the Poisson bracket of
$\cal C^{\infty }(\Omega G) $ by
$$
\lbrace f, g \rbrace (\gamma ) = \omega (\gamma )
(X_f, X_g)
$$
for $f, g\in \cal C^{\infty } (\Omega G)$

(a).The Poisson algebra $Span_{\cal C}\lbrace J_n^a, 1, E\rbrace $
is isomorphic to the Kac-Moody algebra
$\widehat {L\cal G_{\cal C}} = L\cal G_{\cal C} \oplus
\cal C \oplus \cal C.d$, that means,
$$
\begin{array}{ll}
\lbrace J_n^a, J_m^b\rbrace & = f_c^{ab} J_{n+m}^c - im\delta _{n, -m}
\delta ^{ab} \\
\lbrace E, J_n^a\rbrace & = -inJ_n^a
\end{array}
$$

(b). $l_n = \sum _{a,m}J_{n-m}^aJ_m^a,\qquad \lbrace l_n, l_m\rbrace =
(n - m)l_{n+m}$

{\bf Proof:} (1)  Let $\xi (\theta) \in T_{\gamma (\theta )} (\Omega G)$,
then one can choose
$\gamma _t (\theta )$, for small t, depending on t, such that
$
\frac {\partial } {\partial t}\mid _{t=0} \gamma _t (\theta ) = \xi (\theta )
$, then
$$
\begin{array}{ll}
dJ_n^a(\gamma )(\xi (\theta )) & = \xi (\theta ).J_n^a (\gamma )\\
& = \frac {\partial }
{\partial t}\mid _{t=0} J_n^a (\gamma _t (\theta )) \\
& = \frac {\partial } {\partial t}\mid _{t=0}
(\frac {1}{4\pi } \int _0^{2\pi }
< T^a e^{in\theta }
, \gamma _t^{-1} (\theta ) \gamma _t' (\theta )> d\theta ) \\
& = \frac {1} {2\pi } \int _0^{2\pi } \xi (\theta ). <T_ae^{in\theta }
, \gamma '\gamma (\theta )^{-1} > d\theta \\
& = \frac {1} {2\pi } \int _0^{2\pi } <Ad \gamma (\theta )^{-1}.
T_ae^{in\theta },
(\gamma ^{-1} \xi (\theta )) '>d \theta \\
& = \omega (X_{J_n^a} (\gamma ), \xi (\theta ))
\end{array}
$$
$$
\begin{array}{ll}
\xi (\theta ). l_n(\gamma ) & = \frac {\partial }
{\partial t}\mid _{t=0} l_n (\gamma _t (\theta )) \\
& = \frac {\partial } {\partial t}\mid _{t=0}
(\frac {1}{4\pi } \int _0^{2\pi }
<\gamma _t^{-1} (\theta ) \gamma _t' (\theta )e^{in\theta }
, \gamma _t^{-1} (\theta ) \gamma _t' (\theta )> d\theta ) \\
& = \frac {1}{2\pi } \int _0^{2\pi }
<\gamma ^{-1} (\theta ) \gamma '(\theta )e^{in\theta }
, (\gamma ^{-1} (\theta ) \xi (\theta ))'> d\theta  \\
& = \omega (\gamma ^{-1} (\theta ) \gamma '(\theta )e^{in\theta }
, \gamma ^{-1} (\theta ) \xi (\theta )) \\
& = \omega (\gamma )(\gamma '(\theta )e^{in\theta }, \xi (\theta ))
\end{array}
$$
Therefore:
$$
\begin{array}{ll}
X_{J_n^a} (\gamma ) &= -l_{\gamma (\theta) *} Ad \gamma (\theta )^{-1}.
T_ae^{in\theta } \\
& = -R_{\gamma (\theta) *}T_ae^{in\theta }\\
X_{l_n} (\gamma ) &= -e^{in\theta }\gamma '(\theta )
\end{array}
$$

(2) (a) By the definition of the Poisson Bracket on $\Omega G$,
$$
\begin{array}{l}
\quad \lbrace J_n^a, J_m^b\rbrace (\gamma ) \\
 = \omega (\gamma )
(X_{J_n^a}, X_{J_m^b})\\
 =- \frac {1}{2\pi }
\int _0^{2\pi }
< (Ad \gamma ^{-1}(\theta ).e^{in\theta }T^a)'
,Ad \gamma ^{-1}(\theta ).e^{im\theta }T^b>
d\theta \\
 = -\frac {1}{2\pi }
\int _0^{2\pi }
< Ad \gamma ^{-1}(\theta ).(ine^{in\theta }T^a)\\
\qquad  + [ Ad_{\gamma ^{-1}(\theta ).}e^{in\theta }T^a
     , \gamma ^{-1}(\theta )\gamma '(\theta )]
,Ad \gamma ^{-1}(\theta ).e^{im\theta }T^b>
d\theta  \\
 = -\frac {1}{2\pi }
\int _0^{2\pi }
< ine^{in\theta }T^a
  ,e^{im\theta }T^b>
d\theta \\
\quad +
\frac {1}{2\pi }
\int _0^{2\pi }
<  \gamma ^{-1}(\theta )\gamma '(\theta )
 , [ Ad \gamma ^{-1}(\theta ).e^{in\theta }T^a
    , Ad \gamma ^{-1}(\theta ).e^{im\theta }T^b] >
d\theta \\
 =- in\delta _{n, -m}\delta^{ab}
+
\frac {1}{2\pi }
\int _0^{2\pi }
<  \gamma ^{-1}(\theta )\gamma '(\theta )
 , Ad \gamma _{-1} (\theta) [ e^{in\theta }T^a
    ,e^{im\theta }T^b] >
d\theta \\
 = f_c^{ab}
\frac {1}{2\pi }
\int_0^{2\pi }
<  e^{i(m+n)\theta }T^c
  , \gamma '(\theta )\gamma ^{-1}(\theta ) >
d\theta
-
in\delta _{n, -m}\delta ^{ab}\\
 = f_c^{ab}J_{m+n}^c
-
in\delta _{n, -m}\delta ^{ab}
\end{array}
$$
$$
\begin{array}{ll}
\lbrace E, J_n^a\rbrace (\gamma )
& =- \frac {1}{2\pi }
\int _0^{2\pi }
< (Ad \gamma ^{-1}(\theta ).e^{in\theta }T^a)'
 , \gamma ^{-1}(\theta )\gamma '(\theta ) >
d\theta \quad \\
& = -\frac {in}{2\pi }
\int _0^{2\pi }
<e^{in\theta }T^a
 , \gamma '(\theta )\gamma ^{-1}(\theta ) >
d\theta \\
& = -inJ_n^a
\end{array}
$$
the proof of (b) is similar.

    Using this theorem, the symmetry of the Kac-Moody algebra and the conformal
algebra can be seen as two Poisson algebras which consist of
the classical observable functions. It also implies
that after the geometrical quantization, these observable functions will
become operators which give rise to the representation of the Kac-Moody
algebra and the Virasora algebra.

\vskip.3in
\centerline{\bf \S 3. Geometrical Quantization of Wess-Zumino-Witten Model }
\vskip.2in

     Geometrical quantization is a method used to quantify a symplectic
manifold $M$ with an integral symplectic form $\omega$ which will be the
curvature of a line bundle $L$ over $M$. To quantify $M$ is to pick a
complex K$\ddot a$hler structure such that $L$ is a holomorphic line-bundle
and the quantum  Hilbert space is then taken to be the space of
holomorphic sections of $L$. It is well known that if $M$ is taken to be a
coadjoint orbit through an integral weight of a compact Lie group, the
corresponding quantum Hilbert space is precisely the representation space
of a Lie group as described in the Borel-Weil theorem [12,13].

     From the proceeding sections, we have paved the way for the geometrical
quantization. In this section, we will accomplish this procedure to give the
representation of a Loop group which has a complete representation theory
that may be found in [16].

     We reiterate the main results in section 2. The phase space of the WZW
model is a K$\ddot a$hler manifold $\Omega G$ with the integral
symplectic form $\omega$. Moreover, for $\xi \in L\cal G$, we have defined
a $\cal C^{\infty }$-function on $\Omega G$ by $\delta _{\xi }(\gamma )
=\frac {1}{2\pi }\int _0^{2\pi }<\xi ,\gamma '\gamma ^{-1}>d\theta $
with its Hamiltonian vector field $-R_{\gamma *}(\xi )$.
 We take the central extension $\widetilde {\Omega G}$
  of the Loop group $\Omega G=LG/G$ as a $U(1)$-principl
bundle on $\Omega G$. The geometrical quantization line bundle will be
a certain associated bundle of a $U(1)$-representation which is defined
by a positive integer called level.

     First we give a description of the needed connection and its curvature
on the $U(1)$-principal bundle. For $\gamma \in \Omega G$, let $\tilde
\gamma \in
\Omega G =LG/G$ be the lifting of $\gamma $ on $\widetilde {\Omega G}$.
 The connection is defined by the horizontal lifting of the tangent vector
 field on $\Omega G$. Specifically,
for $l_{\gamma *} \xi $, which belongs to the tangent vector at
 $\gamma $ with  $\xi \in \Omega \cal G $, we define
 $ l_{\tilde \gamma *}(\xi ) $ to be the horizontal
lifting of $l_{\gamma *}\xi $ to the point $\tilde \gamma $.

{\bf Theorem 3}:  The curvature form of the above connection is exactly
the symplectic form of $\Omega G$, as needed
in the geometrical quantization. Moreover,
this curvature and connection are left-invariant under the left
multiplication.

{\bf Proof:}
 Assume $l_{\gamma *}\xi , l_{\gamma *}\eta $ to be the two left-invariant
vector field on $\Omega G$, lift them to the point $\tilde \gamma $
by definition. Then the two horizontal lifting vector fields are
$l_{\tilde \gamma *}\xi ,l_{\tilde \gamma *}\eta $
respectively, where $\tilde \gamma $ is the lifting of $\gamma$.
The curvature form is calculated by the following formula:
$$
\begin{array}{ll}
\Omega _{\tilde \gamma }(l_{\tilde \gamma *}\xi ,l_{\tilde \gamma *}\eta )
&=[l_{\tilde \gamma *}\xi ,l_{\tilde \gamma *}\eta ] -l_{\tilde \gamma *}
[\xi ,\eta ]_{\Omega \cal G}\\
&= l_{\tilde \gamma *} [\xi ,\eta ]_{\widetilde {\Omega \cal G}}
-l_{\tilde \gamma *}[\xi ,\eta ]_{\Omega \cal G} \\
&=l_{\tilde \gamma *}\omega (\xi ,\eta )
\end{array}
$$
where $l_{\tilde \gamma *} $ is an $\cal R$-valued vertical tangent vector
field
at $\tilde \gamma $.

       The left-invariant is obvious from the calculation.

	  Let $L_{k}$ be the associated $U(1)$-bundle on $\Omega G$ by the
representation
$$
\rho :U(1) \longrightarrow \cal C \backslash
\lbrace 0 \rbrace
$$
$$
 \rho (e^{i\theta})=e^{ik\theta}
 $$
where k is a positive integer.
Then $L_{k}$ has a natural complex structure such that $L_{k}$ is a
holomorphic line bundle over $(\Omega G, k \omega _{L}, J_{L})$, where
$\omega _{L}, J_{L}$ are the left-invariant  symplectic structure and
the left-invariant complex
structure respectively. From the above theorem, $L_{1}$
is the holomorphic line bundle on $\Omega G$ which will give rise to the
geometrical quantization of $(\Omega G, \omega _{L})$. If we identify
the smooth section of $L_{1}$ with a certain smooth function on
$\widetilde {LG}$
which satisfies the following conditions:
 $$
 \begin{array}{ll}
 (a).\qquad \qquad f(\tilde \gamma g ) &= f(\tilde \gamma ) \qquad for \quad
  g \in G ,
 \tilde \gamma \in \widetilde {LG}\\
 (b). \qquad \qquad f(\tilde \gamma e^{i\theta})&=e^{-i\theta }f(\tilde \gamma
)
 \end{array}
$$
	then we can calculate  the infinitesimal form of the left multiplication
of $ \widetilde {LG}$. This  will give realization of quantum operators which
is the same result as obtained by the geometrical procedure.

{\bf Theorem 4}:   (1) The infinitesimal form of $\widetilde  {LG}$
action on $L_{1}$ gives
rise to the prequantum operator in the geometrical quantization, that is,
for $ \xi \in L\cal G $,
$$
\hat \delta _\xi = \bigtriangledown _{X_{\delta _\xi}} +  i \delta _\xi
$$
    is the infinitesimal representation,
where $\hat \delta _\xi $ is the prequantum operator,
$\bigtriangledown $ is the associated
connection, and $X_{\delta _\xi }$ is the Hamiltonian vector
field of $\delta _\xi $.

             (2) $\hat \delta _\xi $  preserves
   the K$\ddot a$hler polarization.

Therefore, we obtained the representation of $\widetilde {LG} $
 (central extension of the Loop group) and its Lie algebra by
 the procedure of geometrical quantization.

{\bf Proof:} Form the formula,
$Ad\gamma .(\xi ,\lambda ) =
( Ad\gamma .\xi
  , \lambda - \frac {1}{2\pi }
  \int _0^{2\pi }
  <\gamma ^{-1}\gamma '
  ,\xi >
  d\theta
)
$
we obtain that
$$
Ad_{\tilde \gamma ^{-1}} .\widetilde {exp(t\xi )}
= \widetilde {exp(tAd\gamma ^{-1}.\xi )}
e^
{\frac {it}{2\pi }
  \int _0^{2\pi }
    <\gamma '\gamma^{-1}
	, \xi >
  d\theta
}
$$
for t which is small enough, where $\xi \in L\cal G $.

Under the identification of sections and functions satisfying the conditions
described above, $\tilde \gamma $ is a certain lifting of $\gamma $, the
infinitesimal form of $\widetilde {LG}$-action is given by
$$
\begin{array}{ll}
(\xi .f) (\tilde \gamma )
& =
\frac {d} {dt} \mid _{t=0}
f (
  \widetilde {exp (-t\xi )}
  \tilde \gamma
)\\
& =
\frac {d} {dt} \mid _{t=0}
f (
  \tilde \gamma
  Ad_{\tilde \gamma ^{-1}.}
  \widetilde {exp (-t\xi )}
)\\
& =
\frac {d} {dt} \mid _{t=0}
f (
  \tilde \gamma
  \widetilde {exp(-tAd\gamma ^{-1}.\xi )}
  e^
  {\frac {-it}{2\pi }
    \int _0^{2\pi }
      <\gamma '\gamma ^{-1}
	  , \xi >
    d\theta
  }
)\\
& =
\frac {d} {dt} \mid _{t=0}
(
f (
  \tilde \gamma
  \widetilde {exp(-tAd\gamma ^{-1}.\xi )}
)
  e^
  {\frac {it}{2\pi }
    \int _0^{2\pi }
      <\gamma '\gamma ^{-1}
	  , \xi >
    d\theta
  }
)\\
& =
(
\frac {i} {2\pi }
\int _0^{2\pi }
  <\gamma '\gamma ^{-1}, \xi  >
d\theta
)
f(\tilde \gamma )
+
l_{\tilde \gamma *}
(-Ad_{\gamma ^{-1}*}\xi ,) f(\tilde \gamma )\\
& =
(i\delta _{\xi } f)(\tilde \gamma )
+
l_{\tilde \gamma *}
(-Ad_{\gamma ^{-1}*}\xi ,) f(\tilde \gamma )\\
\end{array}
$$

By the definition of connection of the $U(1)$-bundle, $l_{\tilde \gamma *}
(-Ad_{\gamma ^{-1}*}\xi)$ is exactly the horizontal lifting of the
Hamiltonian vector field $X_{\delta _\xi } = - R_{\gamma *}\xi$ to the point
$\tilde \gamma$ . So the proof of (1) above is complete.

{}From the above proof, we know $(\widehat {\delta _\xi }.f)(\tilde \gamma )
=(R_{\tilde \gamma *}\xi )f(\tilde \gamma )$.
the polarization condition is left invariant, given by type(0, 1) vector
fields $L_{\gamma *}(e^{in\theta }T^a)(n < 0)$
 (see the definition of $J_L$ ). The proof is derived from the fact that the
left invariant vector fields and the right invariant vector fields are
commutative, so that the holomorphic sections are preserved by $\widehat
{\delta _\xi }$, and the projective representation of the Loop group is thus
 obtained.

\vskip.4in

\centerline {\bf \S 4.  Moment Map of Gauge Field Theory }
\vskip.in

       In section 1, we have derived the symplectic form on the phase space
of Wess-Zumino-Witten model from the standard symplectic form of the gauge
 field theory on Riemann surfaces. In the following, we will give the moment
 map of gauge field theory and its relationships with the Kohno connection
and the Knizhnik-Zamolodchikov equation in WZW model. This idea was inspired
by reading M. Atiyah's book [2].

   Firstly, we give the general description of gauge field theory on Riemann
surfaces. Let $\Sigma $ be a Riemann surface and $G$ be a compact Lie group,
denote $\cal A$ as the infinite dimensional affine space, composed
by $G$-connections on the trivial $G$ bundle over $\Sigma $, i.e.,
$$
\cal A =\lbrace \alpha \hbox { is  a  $\cal G$-valued  1-form
 on $\Sigma $ } \rbrace
$$

   It is well known that there is a natural symplectic form on $\cal A$,
denoted by $\omega $ , which is given by the following formula:
$$
\omega (\alpha , \beta )= -\frac {1}{2\pi } \int _{\Sigma }Tr (\alpha \wedge
\beta )
$$
where $-Tr$ is the $G$-invariant inner product on $\cal G $, which  also
lead to the invariant property  under the gauge transformation actions.
 In this case, the gauge transformation group is $ Map (\Sigma , G)$,
the smooth maps from $\Sigma $ to $G$, whose Lie algebra is given by
$Map (\Sigma ,\cal G) $, the smooth maps from $\Sigma $ to $\cal G$.
Using the inner product on $\cal G$, the dual space of $Map(\Sigma ,\cal G )$
can be identified with the smooth $\cal G$-valued 2-forms on $\Sigma $,
by integrating over $\Sigma $, while $Map(\Sigma , \cal G )$ is taken as the
smooth $\cal G $-valued 0-form. Namely, we have $Map(\Sigma , \cal G )^* \cong
\bigwedge ^2(\Sigma , \cal G) $.

    Secondly, we calculate the moment map for the action of $Map (\Sigma , G)$
on $\cal A$, as in the M. Atiyah and R. Bott's paper [3].

     A moment map for the action of $Map(\Sigma , G)$ on $\cal A$ is a map
$\mu : \cal A \longrightarrow Map(\Sigma , \cal G )^*$ such that for $\xi \in
Map (\Sigma , \cal G ), v \in T_\alpha \cal A =\cal A , \alpha \in
\cal A$,
$$
<\xi , d\mu _\alpha (v)> = \omega (\xi _\alpha ,v)
$$
where $d\mu _\alpha : T_\alpha \cal A \longrightarrow Map(\Sigma , \cal G)^*$
is the differential map of $\mu $ at $\alpha $, $\xi _\alpha \in T_\alpha
\cal A$ is the vector field defined by $ \xi \in Map(\Sigma , \cal G)$ through
the action of $Map(\Sigma , G)$, and $<,>$ denotes the dual pairing $\cal G$
and $\cal G ^*$.

    Moment map play an important role in classical and quantum  mechanics
by  the Marsden-Weinstein reduction. For more details,
see Marsden-Weinstein [14].

{\bf Theorem 5:} The moment map $\mu : \cal A \longrightarrow Map(\Sigma ,
 \cal G )^*$ exists, concretely,

   (1). If $\Sigma $ has no boundary, then
$$
\mu (\alpha )= d\alpha +\alpha \wedge \alpha =F_\alpha
$$
where $F_\alpha \in \bigwedge ^2 (\Sigma , \cal G) =Map(\Sigma, \cal G )^*$ is
 the curvature form of $\alpha $.

   (2). If $\Sigma $ has a boundary $S=\partial \Sigma $, then
$$
\mu (\alpha ) =F_\alpha -\alpha _S
$$
where $\alpha _S$ is the restriction of $\alpha $ to the boundary $S$, i.e.,
$\alpha_S \in \bigwedge ^1(S, \cal G)$ which is seen as an element of
$Map(\Sigma , \cal G)^* $ in the following sense: for $\xi \in
Map(\Sigma , \cal G) $, $\xi _S$ is the restriction of $\xi $ to the boundary
$S$,
$$
<\alpha _S, \xi >=-\frac {1}{2\pi } \int _S Tr (\xi _S\alpha _S)
$$
{\bf Proof: } By the definition of the moment map, for $\xi \in Map(\Sigma ,
\cal G)=\bigwedge ^0(\Sigma , \cal G), \alpha \in \cal A =\bigwedge (\Sigma ,
\cal G), v\in T_\alpha \cal A =\cal A=\bigwedge ^1(\Sigma , \cal G) $, we have
that the  tangent vector at $\alpha $,
 defined by $\xi \in Map(\Sigma , \cal G)$,  is
$$
\begin{array}{l}
\quad \frac {d}{dt}\mid _{t=0}\lbrace \exp (-t\xi )\alpha \exp (t\xi )+
\exp (-t\xi )d\exp (t\xi )\rbrace \\
=d \xi
+[\alpha , \xi ]=d_\alpha \xi
\end{array}
$$
 Therefore,
we need to verify the following equality:
$$
<\xi , d\mu _\alpha (v)> = \omega (\xi _\alpha ,v)
$$
(1). If $\Sigma $ has no boundary, we only check $\mu (\alpha )=F_\alpha $
is the required moment map.
$$
\begin{array}{ll}
<d\mu _\alpha (\alpha ), \xi >&=\frac {d}{dt}\mid _{t=0}
<F_{\alpha +t v }, \xi>\\
&=<d_\alpha v , \xi >\\
&=-\frac {1}{2\pi }\int _{\Sigma }Tr(d_\alpha v \xi)\\
&=-\frac {1}{2\pi }\int _{\Sigma }Tr(\xi d_\alpha v )\\
&=-\frac {1}{2\pi }\int _{\Sigma }Tr( d_\alpha \xi \wedge v )\\
&=\omega (\xi _{\alpha }, v)
\end{array}
$$
where the fourth equality is obtained by the partial integrating formula for
the oprater $d_\alpha $, $\frac {d}{dt}\mid_{t=0}F_{\alpha +t v}
=\frac {d}{dt}\mid _{t=0}(d(\alpha +tv) +(\alpha +tv)\wedge (\alpha +tv))
=dv +[\alpha , v] =d_\alpha v$.

(2). If $\Sigma $ has a boundary $S=\partial \Sigma $, then
$$
\begin{array}{ll}
<d\mu _{\alpha }(v), \xi >&=\frac {d}{dt}\mid _{t=0}<F_{\alpha +tv}-(\alpha +
tv)_S, \xi >\\
&=<d_\alpha v- v_S , \xi >\\
&=-\frac {1}{2\pi }\int _{\Sigma }Tr(\xi d_\alpha v)+\frac {1}{2\pi }\int _S
Tr(\xi _Sv_S)\\
&=-\frac {1}{2\pi }\int _\Sigma Tr(d_\alpha \xi \wedge v)-\frac {1}{2\pi}
\int _S Tr(\xi _S\wedge v_S) +\frac {1}{2\pi }\int _S Tr(\xi _S\wedge v_S)\\
&=\frac {1}{2\pi}\int _\Sigma Tr(d_\alpha \xi \wedge v)\\
&=\omega (\xi _\alpha , v)
\end{array}
$$

   The most interesting thing is $\Sigma $ with $n$ fixed punctures
$\lbrace P_1 ,\cdots , P_n\rbrace $, at each point $P_i$, equipped with
an appropriate representation $\lambda _i$ of the compact Lie group $G$.
At each point, we also have an evaluation map:
$$
e_{p_i}: \quad Map (\Sigma , G) \longrightarrow G
$$
$$
e_{p_i}(\gamma ) =\gamma (P_i)
$$
where $\gamma \in  Map (\Sigma , G)$, differetial this  map and dual it ,
we obtain an embedding map:
$$
\delta _{P_i}:\quad \cal G^* \longrightarrow Map (\Sigma , \cal G)^*
$$

   By the Borel-Weil theorem, each $\lambda _i $ is corresponding  to an
integal co-adjoint orbit of $G$ in $\cal G ^*$, denoted by $M_i $.
We will give two different Gauss laws of gauge field in the case of $\Sigma $
having no boundary and having  n fixed punctures $\lbrace P_1 ,\cdots , P_n
\rbrace $.  For the former case, the Gauss law is
$$
F_\alpha =0
$$
which is just the Marsden-Weinstein reduction of the moment map $\mu (\alpha )
=F_\alpha $. In the latter case,  associating  each point $P_i$
with a representation
$\lambda _i $ (dominant weight) and an integal co-adjoint orbit $M_i \subset
\cal G ^*$
of $G$, then the Marsden-Weinstein reduction becomes the generalized symplectic
 quotient
$$
\cal M _{\lbrace P_i \rbrace, \lbrace \lambda _i \rbrace }=
\mu ^{-1} [(M_1 + \cdots + M_n)/k ]// Map(\Sigma , G)
$$
where the $ ``//" $ means module the gauge transformation of $Map(\Sigma , G)$.

    For such a gauge field $\alpha \in
\cal M _{\lbrace P_i \rbrace, \lbrace \lambda _i \rbrace }$  its
Gauss  contraints is given by the following:
$$
F_{\alpha }=-\frac {2\pi }{k}\sum _{i, a}\xi _i^a\delta (P_i)T_a
$$
where $F_\alpha  $ is the curvature of $\alpha $, $\lbrace T_a \rbrace $ is an
orthogonal basis of $\cal G$ under the Cartan-Killing form $Tr$,
$\xi _i^a\epsilon _a$ is an element of $M_i$ , $\lbrace \epsilon _a \rbrace $
is the dual basis of $\lbrace T_a\rbrace $ in $\cal G ^*$, $\delta (P_i)$ is
the delta function at $P_i$, which can be taken as an element of
$Map(\Sigma ,\cal G)^*$ by $\delta (P_i)(\xi )=\xi (P_i)$ for $\xi \in
Map (\Sigma , \cal G )$ and the Cartan-Killing form. The above
formula makes sense if using the moment map:
$$
<F_\alpha , \xi >=-\frac {1}{2\pi } \int _\Sigma  Tr(\xi F_\alpha )
=\frac {1}{k} \sum _{i, a }\xi_i^a<\xi (P_i), T_a>
 $$

    Using the local coordinate $(z, \overline z )$ on $\Sigma $, $F_\alpha
=\frac 12 F^a_{z\overline z}(z, \overline z)T_adz\wedge d\overline z$,
 hence, the Gauss
contraints can be written as
$$
F^a_{z\overline z}(z, \overline z)=-\frac {4\pi }{k}
\sum _{i,a}\xi _i^a\delta (z,P_i)
 $$

   We notice that $\lbrace \xi _i^a \rbrace $ is a system of functions with
 degree 1 (coordinate function) on $M_i$, after the geomertric quantization,
these classical observable functions turn out to be the quantum operator
$T_a^i$ acting on the representation space $V_i$ of the representation
$\lambda _i $ [12,13].
This lead to the formula written down in  [20] as follows:
$$
F^a_{z\overline z}(z, \overline z)=-\frac {4\pi }{k}\sum_{i,a}\delta
(z, P_i)T_a^i
$$

    We know that the Green function on compact Riemann surface does not exist,
but on the plane,  there is a Green function due to the following formula:
$$
\partial _z\partial _{\overline z }\ln \mid z-P\mid =2\pi \delta (z, P)
$$
{\bf Theorem 6:} The Gauss constraints have a matrix-valued solution:
$$
\left\{\begin{array}{ll}
A^a_z &=\frac {4}{k}\sum _i\partial _z\ln \mid z-P_i \mid T_a^i=\frac 2k
\sum _i \frac {T_a^i}{z-P_i}\\
A_{\overline z}^{a}&=0
\end{array}\right.
$$
where the connection  is defined on the associated bundle of trivial $G$-bundle
on $\Sigma = \cal R^2\backslash \lbrace P_1, \cdots , P_n \rbrace $ by the
tensor representation
$$\lambda _1\otimes \cdots \otimes \lambda _n : G \longrightarrow GL(V_1
\otimes \cdots \otimes V_n)
$$
{\bf Proof}: This is obtained by the direct calculation:
$$
\begin{array}{ll}
F_A&=\frac 12 F_{z\overline z}dz\wedge d\overline z\\
&=dA + A\wedge A\\
&=-\partial _{\overline z} A_z dz \wedge d\overline z \\
&=-\partial _{\overline z}(\frac 4k\sum _i\partial z \ln \mid z-P_i \mid T_a^i)
dz\wedge d\overline z \\
&=-\frac {4}{k}\sum _i \partial _{\overline z}\partial z\ln \mid z-P_i \mid
 T_a^idz\wedge d\overline z \\
&=-\frac {2\pi }{k}\sum _i \delta(z-P_i)T_a^idz\wedge d\overline z
\end{array}
$$

   This connection $A$ has a special meaning  which not noly gave the
solution of the Gauss constraints, but also leads to the Kohno connection
and the Knizhnik-Zamolodchikov equation in the next section.

\vskip .3in
\centerline {\bf \S 5. Kohno Connection and Knizhnik-Zamolodchikov Equation  }
\vskip .2in

   In the section 4, we have solved the Gauss constraints using Green function.
In this section, we will generalize this connection to Kohno connection [6,11],
whose holonomy gives rise to a representation of the braid group and its
parallel sections will be the solutions of Knizhik-Zamolodchikov equation[4].

    We rewrite the connection of section 4 as
$$
A_q(z)dz=\frac {2}{k}\sum _{p\neq q}\Omega _{pq} d\ln (z-z_p)dz
=\frac {2}{k}\sum _{p\neq q}\frac {\Omega _{pq}}{z-z_p}dz
$$
where $\Omega _{pq}=\sum _a T^p_a\otimes T^q_a $  acts on $V_1\otimes \cdots
\otimes V_n $, only on the $p-th$  and $q-th$ factors nontrivially. In fact,
$\Omega _{pq}$ is the Casimir element of the universal eveloping algebra
$\cal U (\cal  G )$ and satisfied the following equations.
$$
\begin{array}{ll}
[\Omega_{pq}, \Omega_{pr}+\Omega_{qr}]=0  &\quad \mbox { $p < q < r$ }\\[2mm]
[\Omega_{pq}+ \Omega_{pr}, \Omega_{qr}]=0  & \quad \mbox { $p < q < r$ }\\[2mm]
[\Omega_{pq}, \Omega_{rs}] =0  &\quad
\hbox { for distinct p, q, r and s }
\end{array}
$$

    Let us consider the following trivial vector bundle over
$$
\Sigma _n=\lbrace (z_1, \cdots , z_n) \in \cal C^n\mid z_p\neq z_q \quad
\mbox { if $p\neq q$ }\rbrace /S_n
$$
with the fibres $V^{\otimes n}$ (a fixed representation  $V$
 of $\lambda $), where $S_n $ acts on $(z_1, \cdots , z_n) $ by the
 usual permutation of n coordinates. The fundamental group $\pi _1(\Sigma _n)$
is the braid group $B_n$, whose generators are $\sigma _1, \cdots ,
\sigma _{n-1} $, satisfying the following relations:
$$
\begin{array}{l}
\sigma _i\sigma _j=\sigma_j\sigma_i  \qquad \mid i-j \mid \geq 2,\\
\sigma_i\sigma_{i+1}\sigma_i=\sigma_{i+1}\sigma_i\sigma_{i+1} \quad
1\leq i\leq n-2
\end{array}
$$

     The Kohno connection is the generalized one described above, whose
connection 1-form is given by
$$
\begin{array}{ll}
\Theta  &=\frac {k}{(k+h)}\sum A_q(z_q)dz_q \\
&=\frac {2}{k+h}
\sum _{1\leq p < q\leq n}\Omega_{pq}d\ln (z_p-z_q) \\
&=\frac {1}{k+h}\sum _{p\neq q} \frac {T^p_a\otimes T^q_a}{z_p-z_q}(dz_p-dz_q)
\end{array}
$$

    The following theorem is the key point in giving out the relations between
the gauge field theory and the conformal  field theory.

{\bf Theorem 7:} (1). $d\Theta + \Theta \wedge \Theta =0$, namely,
 $\Theta $ is a flat connection. Therefore, the holonomy of $\Theta $
gives the representation of $B_n =\pi _1(\Sigma _n) $, called the monodromy
representation of $B_n$, as in  [6,11].

     (2). The equation of parallel sections
$$
D\psi =0
$$
is just the Knizhnik-Zamolodchikov eqation :
$$
\partial _{z_j}\psi (z_1, \cdots , z_n)=\frac {1}{k+h}\sum _{i\neq j}
\frac {T^i_a\otimes T^j_a}{z_i-z_j}\psi (z_1, \cdots ,z_n).
$$
{\bf Proof :} (1).
$$
\begin{array}{l}
\quad  d\Theta +\Theta \wedge \Theta   \\
=\frac {4}{(k+h)^2}[\sum _{1\leq p<q\leq n}\Omega_{pq}d\ln (z_p-z_q),
\sum _{1\leq r<s\leq n}\Omega_{rs}d\ln (z_r-z_s)]
\end{array}
$$
 From $[\Omega_{pq}, \Omega_{rs}]=0$, for distinct $p, q, r$ and
$s$, we obtain that
$$
\begin{array}{l}
\quad [\sum _{1\leq p<q\leq n}\Omega_{pq}d\ln (z_p-z_q),
\sum _{1\leq r<s\leq n}\Omega_{rs}d\ln (z_r-z_s)]\\
=[\sum _{1\leq p<q\leq n}\Omega_{pq}d\ln (z_p-z_q),
\sum _{1\leq p<s\leq n}\Omega_{ps}d\ln (z_p-z_s)]\\
\quad +[\sum _{1\leq p<q\leq n}\Omega_{pq}d\ln (z_p-z_q),
\sum _{1\leq s<p\leq n}\Omega_{s,p}d\ln (z_s-z_p)]\\
\quad + [\sum _{1\leq p<q\leq n}\Omega_{pq}d\ln (z_p-z_q),
\sum _{1\leq r<q\leq n}\Omega_{r,q}d\ln (z_r-z_q)]\\
\quad +[\sum _{1\leq p<q\leq n}\Omega_{pq}d\ln (z_p-z_q),
\sum _{1\leq q<r\leq n}\Omega_{qr}d\ln (z_q-z_r)]\\
=\sum _{p<q<r}\frac {2}{(z_p-z_q)(z_q-z_r)(z_p-z_r)}\\
\quad \lbrace
[\Omega_{pq}, \Omega_{qr}](z_p-z_r)(dz_p\wedge dz_q-dz_p\wedge dz_r+dz_q
\wedge dz_r)\\
\qquad +[\Omega_{pq}, \Omega_{pr}](z_q-z_r)(dz_p\wedge dz_q-dz_p\wedge dz_r+
dz_q\wedge dz_r)\\
\qquad + [\Omega_{pr}, \Omega_{qr}](z_p-z_q)(dz_p
\wedge dz_q-dz_p\wedge dz_r+dz_q\wedge dz_r)\rbrace \\
=\sum _{p<q<r}
\frac {2(dz_p\wedge dz_q-dz_p\wedge dz_r+dz_q
\wedge dz_r)}{(z_p-z_q)(z_q-z_r)(z_p-z_r)}
\lbrace [\Omega_{pq}, \Omega_{qr}](z_p-z_r)\\
\qquad +[\Omega_{pq}, \Omega_{pr}](z_q-z_r)
+[\Omega_{pr}, \Omega_{qr}](z_p-z_q)\rbrace \\
=\sum _{p<q<r}
\frac {2(dz_p\wedge dz_q-dz_p\wedge dz_r+dz_q
\wedge dz_r)}{(z_p-z_q)(z_q-z_r)(z_p-z_r)}
\lbrace [\Omega_{pq}, \Omega_{pr}+\Omega_{qr}](-z_r)\\
\qquad +[\Omega_{pq}+\Omega_{pr}, \Omega_{qr}]z_p
+[\Omega_{pq}+\Omega_{qr}, \Omega_{pr}]z_q\rbrace
\end{array}
$$

     From $[\Omega_{pq}, \Omega_{pr}+\Omega_{qr}]=
[\Omega_{pq}+ \Omega_{pr}, \Omega_{qr}]=0,$ $ p < q < r,$ we prove that
$$
d\Theta +\Theta \wedge \Theta  =0
$$

  Moreover, the holonomy representation of $\Theta $ is the parallel
displacement along a closed curve $C$, depending only on the homotopy class
 of $C$, especially, let $C$ be an element of $\pi _1(\Sigma _n)=B_n$,
$$
\rho (C)=P\exp (-\oint_C \Theta )
$$
is a well-defined representation of the braid group $B_n$ acting on $V^{\otimes
n}$.

  (2). If $\psi $ is a parallel section, then $\psi $ is a
$V^{\otimes n}$-valued function on $\Sigma_n $ and satisfies the following
 equation:
$$
d\psi +\Theta \psi =0
$$
locally, this equation can be rewritten as follows:
$$
\frac {\partial \psi }{\partial z_j}=\frac {1}{k+h}\sum_{i\neq j}
\frac{T^i_a\otimes T^j_a}{z_i-z_j}\psi(z_1, \cdots ,z_n)
$$

  These equations is exactly the Knizhnik-Zamoldchikov equation [10].
In their paper,
they also showed that the solutions exist in the terms of n-point functions
of the two dimentional conformal field theory with gauge symmetry (WZW model).
In our case, they are the parallel sections of the connection $\Theta $.
The integrability condition of these equations is precisely vanishing of the
curvature $d\Theta +\Theta \wedge \Theta =0$.

\vskip .3in
\centerline {\bf \S 6. Relations to  $R$-Matrix and Quantum Group  }
\vskip .2in

       In Kohno's paper, he had proved that the above holonomy representation
 $\rho $ of braid group $B_n$ is equivalent to the monodromy representation of
braid group $B_n$ which was obtained from the solution of quantum Yang-Baxter
equation (QYBE). We will discuss the latter representation briefly in this
 section.

        Firstly, we introduce the $R$-matrix $R: \cal C^n\otimes \cal C^n
\longrightarrow   \cal C^n\otimes \cal C^n $, which  satisfies the following
equation, called quantum Yang-Baxter equation:
$$
R_{12}R_{13}R_{23}=R_{23}R_{13}R_{12}
$$
where $R_{12}, R_{23}, R_{13}: \cal C^n\otimes \cal C^n
\otimes \cal C^n
\longrightarrow   \cal C^n\otimes \cal C^n \cal C^n$, whose actions are
determined by $R$ and their subscripts. This indicates that $R_{ij}$ acts
on $\cal C^n\otimes \cal C^n\otimes \cal C^n $ like $R$ in the $(i,j)$
factors and like identity in the remaining factor.
 The most common solutions is $A_n$-type, expressed as
follows [5,8]:
$$
R=\sum _{i\neq j}e_{ii}\otimes e_{jj} +q^{\frac 12}\sum _i e_{ii}\otimes
e_{ii}+ (q^{\frac 12}-q^{-\frac 12})\sum _{i < j}e_{ij}\otimes e_{ji}
$$
where $e_{ij}$ is the $n\times n $ matrix unit whose single nonvanishing
component is equal to 1 and it is located in position $(i,j)$.

   Let $\check R= PR$, where $P$ is the permutation operator of $V\otimes V$,
 then, $\check R$ satisfies the following two equations:
$$
\begin{array}{l}
(\check R \otimes I)(I\otimes \check R)(\check R \otimes I)
=(I\otimes \check R)(\check R \otimes I)(I\otimes \check R)\\
\check R ^2 = (q^{\frac 12}-q^{-\frac 12})\check R +1
\end{array}
$$
where the first equality is derived from the QYBE satisfied by $R$, the second
equality means that $\check R$ has two eigenvalues $q^{\frac 12},
-q^{-\frac 12}$.

     By the direct calculation, we obtain that
$$
\check R=\sum _{i\neq j}e_{ij}\otimes e_{ji} +q^{\frac 12}\sum _i e_{ii}\otimes
e_{ii}+ (q^{\frac 12}-q^{-\frac 12})\sum _{i > j}e_{ii}\otimes e_{jj}
$$

     To prove this equation, we choose a basis for $\cal C^n$, denoted by
$\lbrace v_s \rbrace $, then we have:
$$
R(v_s\otimes v_t)=\cases{
v_s\otimes v_t, &\hbox { for $s < t$ }\cr
q^{\frac 12} v_s\otimes v_t, &\hbox { for $s=t$ }\cr
v_s\otimes v_t +(q^{\frac 12}-q^{-\frac 12})v_t\otimes v_s
&\hbox { for $s > t$ }\cr}
$$
Therefore, by $P(v_s\otimes v_t)=v_t\otimes v_s$,
$$
\check R (v_s\otimes v_t)=\cases{
 v_t\otimes v_s, &\hbox { for $s < t$ }\cr
q^{\frac 12} v_t\otimes v_, &\hbox { for $s=t$ }\cr
v_t\otimes v_s +(q^{\frac 12}-q^{-\frac 12})v_s\otimes v_t
&\hbox { for $s > t$ }\cr}
$$
and the expression of $\check R$ is obtained.

    Using $\check R$-matrix, the representation of braid group $B_n$
can be given as follows:
$$
\nu (\sigma _i)=-q^{12}(I\otimes \cdots \otimes \check R\otimes \cdots
\otimes I)
$$
where $\check R$ is located at the $(i, i+1)$-factors, so $\nu(\sigma _i)$
acts on $V^{\otimes n} $ only nontrivially in the $(i, i+1)$-factors.
The relations which $\lbrace \sigma _i \rbrace $ satisfy are satisfied by
$\lbrace \nu _i=\nu (\sigma _i)\rbrace $ due to the equation:
$$
(\check R \otimes I)(I\otimes \check R)(\check R \otimes I)
=(I\otimes \check R)(\check R \otimes I)(I\otimes \check R)
$$

    Moreover, $\lbrace \nu _i \rbrace $ also satisfy the following equation:
$$
\nu_i^2 = (1-q)\nu _i +q
$$
together with the other two relations satisfied by
 $\lbrace \sigma _i\rbrace $,
$\lbrace \nu _i \rbrace $ defines the Hecke algebra of type $A_n$, called the
Terperley-Lieb-Jones algebra [8,9] :
$$
\left\{\begin{array}{l}
\nu_i\nu_{i\pm 1}\nu_i=\nu_{i\pm 1}\nu_i\nu_{i\pm 1}\\
\nu_i\nu_j=\nu_j\nu_i \qquad \mid i-j\mid \geq 2\\
\nu_i^2=(1-q)\nu_i +q
\end{array}\right.
$$

     T. Kohno has stated that the representation using $\check R$-matrix
$$q=\exp (-i\frac {2\pi }{k+h})$$
and the representation through the holonomy of the connection $\Theta $
(section 5) are equivalent up to a conjugation by an invertible linear
transformation of $V^{\otimes n}$. He also proved that the holonomy
representation commutes with the coproduct action of quantum group
$SU_q(n)$ on $V^{\otimes n}$. In fact, the above Hecke algebra is the
centralizer of quantum group $SU_q(n)$.

\vskip.25in
{\bf  Acknowledgement }

    I would like to thank Prof. Qian Min for his encouragement
 and many  useful discussions on this work.  I am also grateful to
 Prof. Guo
Maozheng, Prof. Liu Zhangju and Dr. Wang Zhengdong for their help.

\end{document}